\documentclass[12pt]{amsart}
\usepackage{amscd,amssymb}
\usepackage[graph,frame,poly,arc]{xy}  
\usepackage[plainpages,backref,urlcolor=blue]{hyperref}

\theoremstyle{plain}
\newtheorem{thm}[subsection]{Theorem}
\newtheorem{lem}[subsection]{Lemma}
\newtheorem{prop}[subsection]{Proposition}

\theoremstyle{definition}
\newtheorem{rk}[subsection]{Remark}

\numberwithin{equation}{section}
\newcommand{\A}{{\mathcal A}}

\newcommand{\Q}{\mathbb{Q}}

\newcommand{\C}{\mathbb{C}}
\newcommand{\PP}{\mathbb{P}}

\begin{document}

\title [ Monodromy of the exceptional reflection arrangement of type $G_{31}$]
{On the Milnor monodromy of the exceptional reflection arrangement of type $G_{31}$}

\author[Alexandru Dimca]{Alexandru Dimca$^1$}
\address{Universit\'e C\^ ote d'Azur, CNRS, LJAD, France}
\email{dimca@unice.fr}

\author[Gabriel Sticlaru]{Gabriel Sticlaru}
\address{Faculty of Mathematics and Informatics,
Ovidius University,
Bd. Mamaia 124, 900527 Constanta,
Romania}
\email{gabrielsticlaru@yahoo.com }

\thanks{$^1$ Partially supported by Institut Universitaire de France.}



\begin{abstract} Combining recent results by A. M\u acinic, S. Papadima and  R. Popescu with a spectral sequence and computer aided computations, we determine  the monodromy action on $H^1(F,\C)$, where $F$ denotes the Milnor fiber of the hyperplane arrangement associated to the exceptional irreducible complex reflection group $G_{31}$. This completes the description given by the first author of such monodromy operators for all the other irreducible complex reflection groups.

\end{abstract}
 
\maketitle


\section{Introduction} 

Let $G \subset GL_n(\C)$ be a finite complex reflection group and denote by $\A(G)$ the union of all the reflecting hyperplerplanes of $G$, i.e. the hyperplanes in $\C^n$ fixed by some $g \in G$, $g \ne Id$. For a general reference on complex reflection groups, see \cite{LT} and \cite{OT}.
Since $\A(G)$ is a central hyperplane arrangement, it has a defining equation $f=0$ in $\C^n$, where $f$ 
is a homogeneous polynomial of some degree $d$. One can associate to this setting the Milnor fiber $F(G)$ of the arrangement  $\A(G)$. This is a smooth hypersurface in $\C^n$, defined by $f=1$, and it is  endowed with a monodromy morphism $h: F(G) \to F(G)$, given by $h(x_1,...,x_n)=\exp(2\pi i/d)\cdot (x_1,...,x_n)$, see \cite{BDY, BaSe, BDS, CS, DHA, DL1, DL2, L, MP, Se1, Se2, S2} for related results and to get a feeling of the problems in this very active area.
The study of the induced monodromy operator
$$h^1(G): H^1(F(G),\C) \to H^1(F(G),\C)$$
is the object of the papers \cite{MPP} and \cite{Dreflection}, while in the special case of real reflection groups $G$, there are some additional results on the higher degree monodromy operators  
$$h^j(G): H^j(F(G),\C) \to H^j(F(G),\C)$$
where $j>1$, see  \cite{Se1, Se2, DL2}. In general, not only these monodromy operators are not known, but even the Betti numbers $b_j(F(G))$ are known only in a limited number of cases.

The complex reflection groups have been classified by Shephard and Todd, see \cite{ST}, who showed that there is an infinite series $G(m,p,n)$ of such groups, plus 34 exceptional cases.
The exceptional complex reflection groups in this classification are usually denoted by $G_j$, with $4 \leq j \leq 37$.   
Consider the following polynomial of degree 60 in $S=\C[x,y,z,t]$ 
\begin{equation}
\label{e6}
f=xyzt(x^4-y^4)(x^4-z^4)(x^4-t^4)(y^4-z^4)(y^4-t^4)(z^4-t^4)
\end{equation}
$$((x-y)^2-(z+t)^2) ((x-y)^2-(z-t)^2) ((x+y)^2-(z+t)^2)((x+y)^2-(z-t)^2)$$
$$((x-y)^2+(z+t)^2)((x-y)^2+(z-t)^2)((x+y)^2+(z+t)^2)((x+y)^2+(z-t)^2)$$
$$((x-z)^2+(y+t)^2) ((x-z)^2+(y-t)^2)((x+z)^2+(y+t)^2)((x+z)^2+(y-t)^2)$$
$$((x-t)^2+(y+z)^2)((x-t)^2+(y-z)^2)((x+t)^2+(y+z)^2)((x+t)^2+(y-z)^2).$$
Then the reflection hyperplane arrangement $\A(G_{31})$ consists of 60 hyperplanes in $\C^4$ and is given by $f=0$, see \cite{HR}.
In this note we prove the following result.

\begin{thm}
\label{thmmain}
The  monodromy operator
$$h^1:H^1(F(G_{31}),\C) \to H^1(F(G_{31}),\C)$$
for the exceptional complex reflection group $G_{31}$ is the identity. In particular, the first Betti number $b_1(F(G_{31}))$ is $59$.
\end{thm}
The description of the monodromy operator
$h^1$ for all the other complex reflection groups is given in \cite{Dreflection}, using a  method which cannot be applied to the exceptional group $G=G_{31}$ as explained in 
\cite[Remark 6.2]{Dreflection}. The proof for the case $G=G_{31}$ involves a completely different approach, and this explains why it is written down here as a separate note. Indeed, this proof is close in spirit to our paper \cite{DStCompM}, with computer aided computations playing a key role at several stages. Moreover, the rank of the group $G_{31}$ being $n=4$, we have to deal with higher dimensional singularities and hence with a more complicated spectral sequence than in \cite{DStCompM}, where the case of plane curves corresponding to $n=3$ is discussed, see also Remark \ref{rkdim} (i).

In the second section we recall a spectral sequence approach for the computation of the monodromy of the Milnor fibers of homogeneous polynomials with arbitrary singularities introduced in \cite{Dcomp}, and developed in \cite[Chapter 6]{D1}, \cite{DStCompM}, with several key additions in the joint work of Morihiko Saito and the first author, see \cite{DS1}.
For simplicity, we describe the results only in the case $n=4$, the only case needed in the sequel. However, the approach presented here is very general, and can be applied at least to all free hyperplane arrangements to get valuable information on their first Milnor cohomology.
On the other hand, the success of our method is based on the special properties enjoyed by the mixed Hodge structure on the first cohomology group $H^1(F(G),\Q)$, see the proof of Theorem \ref{thmEV}, and hence this simple approach does not give complete results on
the higher cohomology groups $H^j(F(G),\Q)$, where $j>1$.

In the third section we recall basic facts on the exceptional complex reflection group $G_{31}$, in particular the construction of the basic invariants going back to Maschke \cite{Mas} and the construction of a basis of Jacobian syzygies as described by Orlik and Terao in \cite[Appendix B, pp. 280-281 and p. 285]{OT}.

In the fourth section we describe the algorithm used to determine the monodromy operator
$h^1:H^1(F(G_{31}),\C) \to H^1(F(G_{31}),\C)$.  For this we need a careful study of the second cohomology group of the Koszul complex of the partial derivatives $f_x$, $f_y$, $f_z$ and $f_t$ of the polynomial $f$ of degree 60 given in \eqref{e6}. At several points, this study is done by
using the software SINGULAR \cite{Sing}. To get an idea why computations via SINGULAR are necessary, let us remark that each of the sequences of three homogeneous polynomials in $S$ which are shown to be regular sequences in Proposition \ref{prop1} consists of polynomials of degrees
$28, 12$ and $16$ respectively, and having $136,24$ and respectively $45$ monomials occurring with non zero (and, in fact, quite complicated rational) coefficients. And the final step in the proof is Lemma \ref{lemfinal} where a homogeneous system of 19600 linear equations in 1424 indeterminates is shown to admit only the trivial solution.
The corresponding codes are available at \url{http://math.unice.fr/~dimca/singular.html}.

\medskip

We would like to thank the referee for his useful suggestions.

\section{A spectral sequence for the Milnor monodromy} 
Let $S=\C[x,y,z,t]$ and let $f \in S$ be a homogeneous polynomial of degree $d$ without multiple factors. Consider the global Milnor fiber $F$ of the polynomial $f$ defined by $f(x,y,z,t)=1$ in $\C^4$, with monodromy action $h:F \to F$, 
$$h(x,y,z,t)=\exp(2\pi i/d)\cdot (x,y,z,t).$$

Let $\Omega^j$ denote the graded $S$-module of (polynomial) differential $j$-forms on $\C^4$, for $0 \leq j \leq 4$. The complex $K^*_f=(\Omega^*, d f \wedge)$ is nothing else but the Koszul complex in $S$ of the partial derivatives $f_x$, $f_y$, $f_z$ and $f_t$ of the polynomial $f$. The general theory says that there is a spectral sequence $E_*(f)$, whose first term $E_1(f)$ is computable in terms of the homogeneous components of the cohomology of the Koszul complex $K^*_f$ and whose limit
$E_{\infty}(f)$ gives us the action of monodromy operator  on the graded pieces of $H^*(F,\C)$ with respect to a certain  pole order filtration $P$, see for details \cite{Dcomp}, \cite[Chapter 6]{D1}, \cite{DS1}.

 More precisely,  for any integer $k \in [1,d]$, there is a spectral sequence 
\begin{equation} 
\label{newspsq}
E_1^{s,t}(f)_k=H^{s+t+1}(K^*_f)_{td+k},
\end{equation} 
converging to
\begin{equation} 
\label{limit}
E_{\infty}^{s,t}(f)_k=Gr_P^sH^{s+t}(F,\C)_{\lambda },
\end{equation} 
where $\lambda=\exp (-2 \pi ik/d)$  and $H^{s+t}(F,\C)_{\lambda }$ denotes the associated eigenspace of the monodromy operator.
Moreover, the differential $d_1$ in this spectral sequence is induced by the exterior differentiation of forms, i.e. $d_1 :[\omega] \mapsto [d (\omega)]$. We have the following.

\begin{thm}
\label{thmEV}
Let $D:f=0$ be a reduced degree $d$ surface in $\PP^3$, and let $$\lambda=\exp (-2\pi i k/d)\ne 1,$$  with $k \in (0,d)$ an integer. 
Then $H^1(F,\C)_{\lambda}=0$  if and only if
$$E_2^{1,0}(f)_k =E_2^{1,0}(f)_{k'} = 0,$$
where $k'=d-k$. 
\end{thm}
\proof
Note that on the cohomology group $H^1(F,\C)_{\ne 1}$, which is by definition the direct sum of eigenspaces of $h^1$ corresponding to eigenvalues $\ne 1$, one has $P^2=0$, $F^1 \subset P^1$
and $F^0 = P^0=H^1(F,\C)_{\ne 1}$, where $F^1$ denotes the Hodge filtration, see \cite{DS1}. Let now $\lambda$ be an eigenvalue of the monodromy operator.
The fact that $H^1(F,\C)_{\ne 1}$ is a pure Hodge structure of weight 1, see  \cite{BDS}, \cite{DP}, implies that either
$H^{1,0}(F,\C)_{\lambda } \ne 0$ or $H^{0,1}(F,\C)_{\lambda } \ne 0$. In the first case we get
$Gr_P^{1}H^{1}(F,\C)_{\lambda }\ne 0$, while in the second case we get $H^{1,0}(F,\C)_{\overline \lambda } \ne 0$ and hence $Gr_P^{1}H^{1}(F,\C)_{\overline \lambda }\ne 0$.
The result follows now from the fact that $H^1(K^*_f)=0$ since $f$ is a reduced polynomial.
Indeed, it is known that the smallest $j$ with $H^j(K^*_f) \ne 0$ is precisely the codimension of the singular locus $D_{sing}$ of $D$ in $\PP^3$, see \cite[Theorems A.2.38 and A.2.48]{Eis}. This implies that $E_1^{s,t}(f)_k=0$ for $s+t=0$, and hence 
$$E_2^{1,0}(f)_k =E_{\infty}^{1,0}(f)_k$$
for any $k \in (0,d)$.
\endproof
To check the condition $E_2^{1,0}(f)_k =0$ in practice, we proceed as follows. Consider a 2-form
 \begin{equation}
\label{e9}
\omega= \sum_{1\leq i<j \leq 4}a_{ij}dx_i \wedge dx_j \in \Omega^2,
\end{equation}
where $a_{ij} \in S_{k-2}$ and we use the convention $x=x_1$, $y=x_2$, $z=x_3$, $t=x_4$. 

We have  the following two obvious results.
 \begin{lem}
\label{dfomega=0}
A differential 2-form $\omega \in \Omega^2$ given by \eqref{e9} satisfies $df \wedge \omega=0$ if and only if the following 4 equations in $S$ hold.
$$
(R_1) \ \ a_{34}f_y-a_{24}f_z+a_{23}f_t=0
$$
$$
(R_2) \ \ a_{34}f_x-a_{14}f_z+a_{13}f_t=0
$$
$$
(R_3) \ \ a_{24}f_x-a_{14}f_y+a_{12}f_t=0
$$
$$
(R_4) \ \ a_{23}f_x-a_{13}f_y+a_{12}f_z=0
$$
\end{lem}

 \begin{lem}
\label{lemd=0}
A differential 2-form $\omega \in \Omega^2$ given by \eqref{e9} satisfies $d \omega=0$ if and only if the following 4 equations in $S$ hold.
$$ (E1) \  \ a_{{23}_t}-a_{{24}_z}+a_{{34}_y}=0.$$
$$ (E2) \  \ a_{{13}_t}-a_{{14}_z}+a_{{34}_x}=0.$$
$$ (E3) \  \ a_{{12}_t}-a_{{14}_y}+a_{{24}_x}=0.$$
$$ (E4) \  \ a_{{12}_z}-a_{{13}_y}+a_{{23}_x}=0.$$
\end{lem}

\section{On the exceptional complex reflection group $G_{31}$} 
The exceptional complex reflection group $G_{31}$ has order 46080, and rank 4, hence it acts naturally on $\C^4$. Its basic invariants $f_1$, $f_2$, $f_3$ and $f_4$, of degree $8,12,20$ and $24$ respectively, can be constructed as follows, see \cite{Mas}.
Let 
\begin{equation}
\label{e1}
a=x^4+y^4+z^4+t^4,  \ b=x^2y^2+z^2t^2, \ c= x^2z^2+y^2t^2, \ d=x^2t^2+y^2z^2, \ e=xyzt. 
\end{equation}
We define
\begin{equation}
\label{e2}
A_1=a+6(-b-c-d), \ A_2= a+6(-b+c+d),  \  A_3= a+6(b-c+d), 
\end{equation}
$$ \ A_4=a+6(b+c-d), \ A_5=-2a-24e, \ A_6=-2a+24e.$$
Let $s_i$ be the $i$-th symmetric function in the variables $A_i$, such that
$$\prod_{i=1,6} (u+A_i)=\sum _{j=0,6}s_ju^{6-j},$$
in the polynomial ring $S[u]$, with $s_0=1$. Define the following homogeneous polynomials in $S$
\begin{equation}
\label{e3}
F_8=-\frac{1}{6}s_2, \ F_{12}=-\frac{1}{4}s_3, \ F_{20}=\frac{1}{12}s_5,
\end{equation}
and note that $s_4=9F_8^2$. Define new polynomials as in \cite{OT}, p. 285, by the following
\begin{equation}
\label{e4}
f_1=F_8, \ f_2=F_{12}, \ f_3=  F_{20}, \ f_4=\frac{1}{265531392}Hess(f_1),
\end{equation}
where $Hess(f_1)$ is the determinant of the Hessian matrix $H(f_1)$ of $f_1$.

 Moreover, as any reflection hyperplane arrangement, $\A(G_{31})$ is a free arrangement, see \cite[Theorem 6.60]{OT}, \cite{Te}. This means that the graded $S$-module of Jacobian syzygies for $f$, namely
$$AR(f)=\{ r=(r_1,r_2,r_3,r_4) \in S^4 \ : \ r_1f_x+r_2f_y+r_3f_z+r_4f_t=0 \},$$
which is the same as the $S$-module of derivations killing $f$, is free of rank $3$. A basis for this module can be computed starting from the basic invariants $f_1,f_2, f_3,f_4$ as follows, see \cite[Theorem 6.53 and Appendix B, pp. 280-281 and p. 285]{OT}. Let $A(f_1)$ be the adjoint of the Hessian matrix  $H(f_1)$, namely a matrix such that
$$A(f_1)H(f_1)=H(f_1)A(f_1)=Hess(f_1)I_4,$$
with $I_4$ the unit matrix of order 4. This matrix is, up to the factor $Hess(f_1)$, the inverse matrix of $H(f_1)$.
Let $J(f)$ be the $4 \times 4$ matrix having on the $j$-th column the first order partial derivatives of $f_j$ with respect to $x,y,z$ and $t$. Define a new $4 \times 4$ matrix $C$, as described in  \cite[Appendix B, p. 284]{OT} for the last six exceptional groups. For the group $G_{31}$, the matrix $C$ has the first, the third and the fourth rows the same as the unit matrix $I_4$, and the second row is given by 
\begin{equation}
\label{e5}
0, \ f_4, \ -\frac{1}{5}f_1 , -\frac{1}{1620}f_2,
\end{equation}
see \cite[Appendix B, p. 285]{OT}.
Then we set $B=A(f_1)J(f)C$ and it follows from \cite[Theorem 6.53]{OT} that the first column in $B$ is a multiple of the Euler derivation. Let $b_{11}$ be the entry of this matrix $B$ situated on the first row and the first column, set $g=b_{11}/x$ and let $D=g^{-1}B$, a new matrix with entries in the polynomial ring $S$. 
The first column $D[1]$ in the matrix $D$, which is a normalized version of $B$, is now given by $x,y,z,t$, and the corresponding derivation of the polynomial ring $S$, denoted again by $D[1]$,  is precisely the Euler derivation.
In particular, one has $D[1]f=60f$. Note that the other columns $D[j]$  of the matrix $D$, for $j=2,3,4$, satisfy similar relations $D[j]f=g_jf$ for some homogeneous polynomials $g_j \in S$, for $j=2,3,4$.
Define a new matrix $E$, whose columns $E[j]$ are given by the following relations.
 \begin{equation}
\label{e7}
E[j]=D[j]-\frac{g_j}{60}D[1],
\end{equation}
for $j=1,2,3,4$, where we set $g_1=0$. It follows that $E[j]f=0$ for $j=2,3,4$, hence we have obtained 3 syzygies for $f$. Moreover, a direct computation using Singular shows that
 \begin{equation}
\label{e8}
\det E= -486 f
\end{equation}
and hence, in view of Saito's Criterion \cite[Theorem 4.19]{OT}, $f$ is free with exponents $1,13,17,29$. In other words, the columns $E[j]$ for $j=2,3,4$ give a basis for the free $S$-module $AR(f)$. Note that $\deg E[2]=29$, $\deg E[3]=13$ and $\deg E[4]=17$. The set of these degrees, to which we must add 1, the degree of the Euler derivation, correspond to the coexponents of $G_{31}$ as listed in \cite[Table B.1., p. 287]{OT}. The entries of the matix $E$ are too complicated to display here, see Remark \ref{rkcoeffsyz}.

Finally we discuss the monodromy operator $h^1:H^1(F(G_{31}),\C) \to H^1(F(G_{31}),\C)$. To compute this operator, we can take a general hyperplane section $V(G_{31}): g=0$ of the projective surface $D:f=0$.
Indeed, the first Milnor monodromy for $g$ and $f$ coincide, see \cite[Theorem 4.1.24]{D1}. Table C.12 in \cite{OT} shows that the corresponding curve $V(G_{31})$ in $\PP^2$ has degree $d=60$, 360 double points, 320 triple points and 30 points of multiplicity 6, corresponding to the isotropy group $G(4,2,2)$. It follows that $h^1$ can have only eigenvalues of order $1,2,3$ and $6$, since any eigenvalue of $h^1$ has to be a root of at least one of the local Alexander polynomials associated to the singularities of $V(G_{31})$, see \cite{L1}, \cite[Corollary 6.3.29]{D1}. The eigenvalues of order $2$ and $3$ are excluded by the results in \cite{MPP}.

To prove Theorem \ref{thmmain}, it remains to exclude the eigenvalues of order 6, and to do this, we apply  Theorem \ref{thmEV}, for $d=60$ and $k=10$ or $k=50$.
Note that for $k=10$ we have 
$E^{1,0}_1(f)_{10}=H^2(K^*_f)_{10}=0.$
Indeed, by the above computations one has 
$AR(f)_8=0$ as all the generating syzygies $E[j]$ have degree $\geq 13$, and then Lemma \ref{dfomega=0}  implies that $H^2(K^*_f)_{10}=0$. 
Therefore it remains to show that $E^{1,0}_2(f)_{50}=0$, and this is done in the next section.
 \begin{rk}
\label{rkdim}
(i) In principle one can try to compute $h^1$ by applying the approach described in \cite{DStCompM} to the curve $V(G_{31}): g=0$ defined above. However, this curve is no longer free, and it seems to us better from the computational point of view to work with the free surface $D$, having non isolated singularities, than to work with the non free curve $V(G_{31})$ with isolated singularities. This is due to the existence of the nice procedure to compute a basis for $AR(f)$ described above.

\noindent (ii) A computation of $h^1$ using the curve $V(G_{31}): g=0$ defined above can be found in \cite{BDY}. The main techniques there are some subtle vanishing results for the cohomology of perverse sheaves on affine varieties, coupled again with computer aided
computation via Singular needed to show that some surfaces are affine.

\noindent (iii) One can avoid using the results in \cite{MPP} and exclude the eigenvalues of order 2 and 3 by the same approach as we use in the case of eigenvalues of order 6. These computations correspond to $k \in \{20,30,40\}$. The cases  $k \in \{20,30\}$ are trivial, since by an obvious modification of Theorem \ref{thm1} below, one has  $H^2(K^*_f)_k=0$ in such situations. The case $k=40$ can be settled by  following the algorithm described below: with the notation of Theorem \ref{thm1} we have in this case $\deg A'_1=8$ and $A'_2=A'_3=0$, i.e. much simpler computations give the result. The corresponding Singular code $g31cubic$ is available at \url{http://math.unice.fr/~dimca/singular.html}.
\end{rk}

\section{The algorithm}

In this section we assume that the coefficients $a_{ij}$ of the  differential 2-form $\omega \in \Omega^2$ given by \eqref{e9} are homogeneous polynomials in $S$ of degree 48, which corresponds to $\deg \omega=50$.

We denote the columns of the matrix $E$ above as follows:
$E[2]=(m_1,n_1,p_1,q_1)$, $E[3]=(m_2,n_2,p_2,q_2)$, $E[4]=(m_3,n_3,p_3,q_3)$.
All the polynomials $m_i$ are divisible by $x$ (since $f_y$, $f_z$ and $f_t$ are clearly divisible by $x$ and $f_x$ is not) and we set $m_i=xm'_i$ for $i=1,2,3$.
Similarly one has $n_i=yn'_i$, $p_i=zp'_i$ and $q_i=tq'_i$ for $i=1,2,3$.
 \begin{prop}
\label{prop1}
Each of the four sequences $(m'_1,m'_2,m'_3)$, $(n'_1,n'_2,n'_3)$, $(p'_1,p'_2,p'_3)$ and $(q'_1,q'_2,q'_3)$ is a regular sequence in $S$.
\end{prop}
\proof
To show that $(m'_1,m'_2,m'_3)$ is a regular sequence, it is enough to check that the dimension of the corresponding zero set in $\C^4$ is the expected one, i.e.
$$\dim V(m'_1,m'_2,m'_3)=1.$$
 And this can be checked by Singular. 
\endproof
 \begin{rk}
\label{rkcoeffsyz}
The four sequences $(m'_1,m'_2,m'_3)$, $(n'_1,n'_2,n'_3)$, $(p'_1,p'_2,p'_3)$ and $(q'_1,q'_2,q'_3)$ of polynomials in $S$ are listed in the print out of a Singular code here.\\
\url{http://math.unice.fr/~dimca/singular.html}
\end{rk}

Note that the relation $R_1$ in Lemma \ref{dfomega=0} is a syzygy in $AR(f)$ having the first component zero. Hence there are polynomials $P,Q,R$ in $S$ such that
$$R_1=PE[2]+QE[3]+RE[4]$$
and $Pm_1+Qm_2+Rm_3=0$. By dividing the last relation by $x$, we get
\begin{equation}
\label{r5}
Pm'_1+Qm'_2+Rm'_3=0.
\end{equation}
By Proposition \ref{prop1}, it follows that $(P,Q,R)$ must be an $S$-linear combination of Koszul relations
for the regular sequence $(m'_1,m'_2,m'_3)$. Consider the following syzygies in $AR(f)$.
\begin{equation}
\label{r6}
T_1=m'_2E[4]-m'_3E[3], \ T_2=m'_1E[4]-m'_3E[2], \ T_3=m'_1E[3]-m'_2E[2].
\end{equation}
More precisely, one has
$$T_1=(0, y(m'_2n'_3-m'_3n'_2),z(m'_2p'_3-m'_3p'_2),t(m'_2q'_3-m'_3q'_2))$$
and similar formulas for $T_2$ and $T_3$. One has $\deg T_1=29$, $\deg T_2=45$, $\deg T_3=41$.
It follows that there are homogeneous polynomials $A_1,A_2,A_3$ in $S$ with $\deg A_1= 19$, $\deg A_2= 3$, $\deg A_3= 7$
such that
\begin{equation}
\label{r7}
R_1=A_1T_1+A_2 T_2 +A_3T_3.
\end{equation}
To express $R_2$ as a combination of simpler syzygies, we have to consider for the same reason as above the syzygies in $AR(f)$ given by
\begin{equation}
\label{r8}
U_1=n'_2E[4]-n'_3E[3], \ U_2=n'_1E[4]-n'_3E[2], \ U_3=n'_1E[3]-n'_2E[2].
\end{equation}
More precisely, one has
$$U_1=(-x(m'_2n'_3-m'_3n'_2),0, z(n'_2p'_3-n'_3p'_2),t(n'_2q'_3-n'_3q'_2))$$
and similar formulas for $U_2$ and $U_3$. One has $\deg U_1=29$, $\deg U_2=45$, $\deg U_3=41$.
It follows that there are homogeneous polynomials $B_1,B_2,B_3$ with $\deg A_j=\deg B_j$, such that
\begin{equation}
\label{r9}
R_2=B_1U_1+B_2 U_2 +B_3U_3.
\end{equation}
However, the coefficient of $f_y$ in the relation $R_1$ coincides with the coefficient of $f_x$ in the relation $R_2$, hence we have
$$yA_1(m'_2n'_3-m'_3n'_2)+yA_2(m'_1n'_3-m'_3n'_1)+yA_3(m'_1n'_2-m'_2n'_1)=$$
$$-xB_1(m'_2n'_3-m'_3n'_2)-xB_2(m'_1n'_3-m'_3n'_1)-xB_3(m'_1n'_2-m'_2n'_1).$$
In other words,
\begin{equation}
\label{r10}
(yA_1+xB_1)(m'_2n'_3-m'_3n'_2)+(yA_2+xB_2)(m'_1n'_3-m'_3n'_1)+(yA_3+xB_3)(m'_1n'_2-m'_2n'_1)=0,
\end{equation}
i.e. a syzygy for the $2$-minors $MN_{ij}$ in the $2 \times 3$ matrix $MN$ having the first row $(m'_1,m'_2,m'_3)$ and the second row $(n'_1,n'_2,n'_3)$.
\begin{rk}
\label{rk1}
The relation $(R_1)$ implies that $a_{34}$ is divisible by $g=xy(x^4-y^4)$, since both $f_z$ and $f_t$ are divisible by $g$, but not $f_x$ or $f_y$ as they have no multiple factors.
A direct computation by Singular shows that all the minors $MN_{ij}$ are also divisible by $g_0=x^4-y^4$.
Therefore, when we want to determine $a_{34}$ as an $S$-linear combination of $MN_{ij}$, we can simplify by $g_0$.
\end{rk}
By the above remark, we can consider new polynomials $MN'_{ij}=(m'_in'_j-m'_jn'_i)/g_0$
for $1 \leq i <j \leq 3$ and hence \eqref{r10} can be written as
\begin{equation}
\label{r11}
(yA_1+xB_1)MN'_{23}+(yA_2+xB_2)MN'_{13} +(yA_3+xB_3)MN'_{12}=0.
\end{equation}
\begin{lem}
\label{lem12}
The minimal degree syzygy  
$r_1MN'_{23}+r_2MN'_{13}+r_3MN'_{12}=0$ has multidegree $(\deg r_1, \deg r_2, \deg r_3)=(24,8,12)$.
\end{lem}
\proof
Direct computation by Singular.
\endproof
It follows that $yA_1+xB_1=yA_2+xB_2=yA_3+xB_3=0$, in other words, there are homogeneous polynomials
$A'_1,A'_2,A'_3$ such that $A_i=xA'_i$ and $B_i=-yA'_i$ for $i=1,2,3$.
Consider now the third syzygy $R_3$, having the third coordinate trivial. We define  new syzygies in $AR(f)$ by 
\begin{equation}
\label{r12}
V_1=p'_2E[4]-p'_3E[3], \ V_2=p'_1E[4]-p'_3E[2], \ V_3=p'_1E[3]-p'_2E[2].
\end{equation}
More precisely, one has
$$V_1=(-x(m'_2p'_3-m'_3p'_2),-y(n'_2p'_3-n'_3p'_2), 0, t(p'_2q'_3-p'_3q'_2))$$
and similar formulas for $V_2$ and $V_3$. One has $\deg V_1=29$, $\deg V_2=45$, $\deg V_3=41$.
It follows that there are homogeneous polynomials $C_1,C_2,C_3$ with $\deg A_j=\deg C_j$, such that
\begin{equation}
\label{r13}
R_3=C_1V_1+C_2 V_2 +C_3V_3.
\end{equation}
However, the coefficient of $f_z$ in $R_1$ coincides with the opposite of the coefficient of $f_x$ in $R_3$, hence we have
$$zA_1(m'_2p'_3-m'_3p'_2)+zA_2(m'_1p'_3-m'_3p'_1)+zA_3(m'_1p'_2-m'_2p'_1)=$$
$$xC_1(m'_2p'_3-m'_3p'_2)+xC_2(m'_1p'_3-m'_3p'_1)+xC_3(m'_1p'_2-m'_2p'_1).$$
In other words,
\begin{equation}
\label{r14}
(zA_1-xC_1)(m'_2p'_3-m'_3p'_2)+(zA_2-xC_2)(m'_1p'_3-m'_3p'_1)+(zA_3-xC_3)(m'_1p'_2-m'_2p'_1)=0,
\end{equation}
i.e. a syzygy for the $2$-minors $MP_{ij}$ in the matrix $MP$ having the first row $(m'_1,m'_2,m'_3)$ and the second row $(p'_1,p'_2,p'_3)$.

\begin{lem}
\label{lem13}
\begin{enumerate}

\item The minors $(MP_{23},MP_{13},MP_{12})$ are all divisible by $x^4-z^4$.

\item Set $MP'_{ij}=MP_{ij}/(x^4-z^4)$. Then the  minimal degree syzygy involving $(MP'_{23},MP'_{13},MP'_{12})$ has multidegree $(24,8,12)$.

\end{enumerate}
\end{lem}
\proof
Direct computation by Singular.
\endproof
It follows as above that one has $C_i=zA'_i$ for $i=1,2,3$.
Consider now the last syzygy $R_4$, having the fourth coordinate trivial. We define the syzygies
\begin{equation}
\label{r20}
W_1=q'_2E[4]-q'_3E[3], \ W_2=q'_1E[4]-q'_3E[2], \ W_3=q'_1E[3]-q'_2E[2].
\end{equation}
More precisely, one has
$$W_1=(-x(m'_2q'_3-m'_3q'_2),-y(n'_2q'_3-n'_3q'_2), -z(p'_2q'_3-p'_3q'_2),0)$$
and similar formulas for $W_2$ and $W_3$. One has $\deg W_1=29$, $\deg W_2=45$, $\deg W_3=41$.
It follows that there are homogeneous polynomials $D_1,D_2,D_3$ with $\deg A_j=\deg D_j$, such that
\begin{equation}
\label{r21}
R_4=D_1W_1+D_2 W_2 +D_3W_3.
\end{equation}
Moreover, the coefficient of $f_t$ in $R_1$ coincides with the coefficient of $f_x$ in $R_3$, hence we have
$$tA_1(m'_2q'_3-m'_3q'_2)+tA_2(m'_1q'_3-m'_3q'_1)+tA_3(m'_1q'_2-m'_2q'_1)=$$
$$-xD_1(m'_2q'_3-m'_3q'_2)-xD_2(m'_1q'_3-m'_3q'_1)-xD_3(m'_1q'_2-m'_2q'_1).$$
In other words,
\begin{equation}
\label{r22}
(tA_1+xD_1)(m'_2q'_3-m'_3q'_2)+(tA_2+xD_2)(m'_1q'_3-m'_3q'_1)+(tA_3+xD_3)(m'_1q'_2-m'_2q'_1)=0,
\end{equation}
i.e. a syzygy for the $2$-minors $MQ_{ij}$ in the matrix $MQ$ having the first row $(m'_1,m'_2,m'_3)$ and the second row $(q'_1,q'_2,q'_3)$.

\begin{lem}
\label{lem23}
\begin{enumerate}

\item The minors $(MQ_{23},MQ_{13},MQ_{12})$ are all divisible by $x^4-t^4$.

\item Set $MQ'_{ij}=MQ_{ij}/(x^4-t^4)$. Then the  minimal degree syzygy involving $(MQ'_{23},MQ'_{13},MQ'_{12})$ has multidegree $(24,8,12)$.

\end{enumerate}
\end{lem}
\proof
Direct computation by Singular.
\endproof
It follows as above that one has $D_i=-tA'_i$ for $i=1,2,3$. This proves the following result.
 \begin{thm}
\label{thm1}
A differential 2-form $\omega$ of degree 50 given by \eqref{e9} satisfies $df \wedge \omega=0$ if and only if there are homogeneous polynomials $A'_1 \in S_{18}$, $A'_2 \in S_2$ and $A'_3 \in S_6$ such that
$$R_1=x(A'_1T_1+A'_2T_2+A'_3T_3), \ R_2=-y(A'_1U_1+A'_2U_2+A'_3U_3),$$
$$R_3=z(A'_1V_1+A'_2V_2+A'_3V_3), \ R_4=-t(A'_1W_1+A'_2W_2+A'_3W_3).$$
\end{thm}
\proof
One should check that any $a_{ij}$ which occurs in two distinct relations $R_k$ and $R_{k'}$
gets the same values by this construction. As an illustration, let's check this for $a_{12}$.
In the syzygy $R_3$, the polynomial $a_{12}$ occurs on the last coordinate, hence its value is
$$zA'_1[t(p'_2q'_3-p'_3q'_2)]+zA'_2[t(p'_1q'_3-p'_3q'_1)]+zA'_3[-t(p'_2q'_1-p'_1q'_2)].$$
In the syzygy $R_4$, the polynomial $a_{12}$ occurs on the third coordinate, hence its value is
$$-tA'_1[-z(p'_2q'_3-p'_3q'_2)]-tA'_2[-z(p'_1q'_3-p'_3q'_1)]-tA'_3[z(p'_2q'_1-p'_1q'_2)],$$
hence exactly the same value as above.
\endproof
The above formulas show that we have the equalities
	\begin{itemize}
		\item[i)] $a_{12}=zt(A'_1PQ'_{23}+A'_2PQ'_{13}+A'_3PQ'_{12}),$
		\item[ii)] $a_{23}=xt(A'_1MQ'_{23}+A'_2MQ'_{13}+A'_3MQ'_{12}),$
		\item[iii)] $a_{13}=-yt(A'_1NQ'_{23}+A'_2NQ'_{13}+A'_3NQ'_{12}).$
	\end{itemize}	

The polynomials $A'_i$ are unique, since the  minimal degree syzygy involving say  $(PQ'_{23},PQ'_{13},PQ'_{12})$ has multidegree $(24,8,12)$ by Lemma \ref{lem23}. It follows that
$$\dim H^2(K^*_f)_{50}=\dim S_{18}+\dim S_2+\dim S_6=1330+10+84=1424.$$
A differential form $\omega \in  H^2(K^*_f)_{50}$ will survive to give an element in $E^{1,0}_2(f)_{50}$ if and only if it satisfies the equations
$(E_i)$ for $i=1,2,3,4$ in Lemma \ref{lemd=0}.
Any such equation $(E_i)$ is in fact a system of $19600=\dim S_{47}$ linear equations in the $1424$  indeterminates given by the coefficients of the polynomials $A'_1$, $A'_2$ and $A'_3$.
Our proof is completed by the following result.
\begin{lem}
\label{lemfinal}
The system of linear equations corresponding to the equation $(E_4)$
in Lemma \ref{lemd=0} has only the trivial solution. In particular one has
$E^{1,0}_2(f)_{50}=0$.
\end{lem}
\proof
Direct, rather lengthy, computation by Singular. The corresponding code $g31$ is available at 
\url{http://math.unice.fr/~dimca/singular.html}
\endproof

\end{document}